\documentclass[a4paper,11pt]{amsart}
\usepackage[colorlinks, linkcolor=blue,anchorcolor=Periwinkle,
    citecolor=blue,urlcolor=Emerald]{hyperref}
\usepackage[all]{xy}
\SelectTips{cm}{}

\usepackage{graphicx}
\usepackage{psfrag}

\usepackage{tikz}
\usepackage{tikz-cd}

\newcommand{\ie}{{\em i.e.}\ }
\newcommand{\cf}{{\em cf.}\ }

\newcommand{\ko}{\: , \;}

\newcommand{\ol}[1]{\overline{#1}}

\numberwithin{equation}{subsection}
\newtheorem{theorem}[subsection]{Theorem}
\newtheorem{classification-theorem}[subsection]{Classification Theorem}
\newtheorem{decomposition-theorem}[subsection]{Decomposition Theorem}
\newtheorem{proposition-definition}[subsection]{Proposition-Definition}
\newtheorem{periodicity-conjecture}[subsection]{Periodicity Conjecture}
\newtheorem{lemma}[subsection]{Lemma}
\newtheorem{keylemma}[subsection]{Key lemma}
\newtheorem{proposition}[subsection]{Proposition}
\newtheorem{corollary}[subsection]{Corollary}

\newtheorem{remark}[subsection]{Remark}
\newtheorem{exercise}[subsection]{Exercise}
\newtheorem{definition}[subsection]{Definition}

\newcommand{\Mod}{\mathrm{Mod}\,}

\newcommand{\proj}{\mathrm{proj}\,}

\renewcommand{\mod}{\mathrm{mod}\,}

\newcommand{\per}{\mathrm{per}\,}
\newcommand{\pvd}{\mathrm{pvd}\,}

\newcommand{\add}{\mathrm{add}\,}

\newcommand{\Hom}{\mathrm{Hom}}
\newcommand{\RHom}{\mathrm{RHom}}
\newcommand{\Ext}{\mathrm{Ext}}

\newcommand{\fib}{\mathrm{fib}\,}
\newcommand{\cof}{\mathrm{cof}\,}

\newcommand{\Z}{\mathbb{Z}}

\newcommand{\la}{\leftarrow}
\newcommand{\iso}{\xrightarrow{_\sim}}

\newcommand{\Si}{\Sigma}
\newcommand{\Ga}{\Gamma}
\newcommand{\La}{\Lambda}

\newcommand{\ten}{\otimes}
\newcommand{\lten}{\overset{\boldmath{L}}{\ten}}

\newcommand{\ca}{{\mathcal A}}
\newcommand{\cb}{{\mathcal B}}
\newcommand{\cc}{{\mathcal C}}
\newcommand{\cC}{{\mathcal C}}
\newcommand{\cd}{{\mathcal D}}

\newcommand{\ch}{{\mathcal H}}

\newcommand{\cm}{{\mathcal M}}

\newcommand{\eps}{\varepsilon}
\renewcommand{\phi}{\varphi}

\renewcommand{\hat}[1]{\widehat{#1}}

\renewcommand{\tilde}[1]{\widetilde{#1}}

\usepackage{tikz}
\usetikzlibrary{cd}
\usetikzlibrary{decorations.pathreplacing}
\usetikzlibrary{matrix,arrows}

\usepackage{tkz-graph}
\GraphInit[vstyle = Shade]
\tikzset{
  LabelStyle/.style = {rectangle, rounded corners, draw,
                        minimum width = 2em, fill = yellow!50,
                        text = blue, font = \Large\bfseries},
  VertexStyle/.append style = { inner sep=5pt,
                                font = \Large\bfseries},
  EdgeStyle/.append style = {->, bend left} }

\begin{document}

\date{August 2021}

\title{An introduction to relative Calabi--Yau structures}

\author{Bernhard Keller}
\address{
Universit\'e Paris Cit\'e, UFR de math\'ematiques, CNRS IMJ--PRG, Sorbonne Universit\'e \\
8 place Aur\'elie Nemours, 75013 Paris, France}
\email{bernhard.keller@imj-prg.fr}
\urladdr{https://webusers.imj-prg.fr/~bernhard.keller/}

\author{Yu Wang}
\address{School of Mathematics and Statistics,
Taiyuan Normal University,
Jinzhong 030619,
PR China
}
\email{dg1621017@smail.nju.edu.cn}
\address{
Universit\'e Paris Cit\'e, UFR de math\'ematiques, CNRS IMJ--PRG, Sorbonne Universit\'e \\
8 place Aur\'elie Nemours, 75013 Paris, France}
\email{yu.wang@imj-prg.fr}






\begin{abstract}
These are notes taken by the second author 
for a series of three lectures by the first author on absolute and relative Calabi--Yau 
completions and Calabi--Yau structures given at the workshop of the International Conference on 
Representations of Algebras which was held online in November 2020. Such structures
are relevant for (higher) representation theory as well as for the categorification
of cluster algebras with coefficients. 
After a quick reminder on dg categories and their Hochschild and cyclic homologies,
we present examples of absolute and relative Calabi--Yau completions
(in the sense of Yeung). In many examples, these are related to higher preprojective algebras
in the sense of Iyama--Oppermann. 
We conclude with the definition of relative (left and right) Calabi--Yau structures after 
Brav--Dyckerhoff.
\end{abstract}

\maketitle

\section{Introduction}

In 1957, the Italo-American geometer Eugenio Calabi conjectured \cite{Calabi57} 
that each K\"ahler manifold whose
first Chern class vanishes admits a Ricci flat metric. His conjecture was proved twenty years later
by Shing-Tung Yau \cite{Yau78}. Algebraic varieties with trivial canonical bundle are now called Calabi--Yau
varieties. The proper ones among these are characterized by the fact that their bounded derived category admits a Serre
functor isomorphic to a power (the dimension) of the suspension functor.
Following Kontsevich, a Hom-finite (algebraic) triangulated category with this property is now called a
Calabi--Yau triangulated category. In these notes, we are concerned with a relative version of this
notion that was first sketched by To\"en \cite{Toen14} in 2014 and fully developed by 
Brav and Dyckerhoff in \cite{BravDyckerhoff19, BravDyckerhoff18}. One of the key
features of their notion of (left) relative Calabi--Yau structure is a gluing construction
analogous to that in cobordism of manifolds. Wai-Kit Yeung showed how to construct
large classes of examples using relative Calabi--Yau completions in \cite{Keller11b, Yeung16}
and advocated the idea that these should be viewed as noncommutative conormal
bundles. This was justified using Kontsevich--Rosenberg's criterion by
Bozec--Calaque--Scherotzke in \cite{BozecCalaqueScherotzke20}.

In the representation theory of quivers and finite-dimensional algebras, the motivations
for studying relative Calabi--Yau structures come from at least three sources:
\begin{itemize}
\item[-] applications in the study of Fukaya categories (cf. for example
Brav--Dyckerhoff's \cite{BravDyckerhoff19, BravDyckerhoff18}),
\item[-] the categorification of cluster algebras with coefficients
(as in the work of Geiss--Leclerc--Schr\"oer \cite{GeissLeclercSchroeer11b}, Leclerc \cite{Leclerc16},
Jensen--King--Su \cite{JensenKingSu16}, Pressland \cite{Pressland17, Pressland17b, Pressland19, Pressland20}  \ldots as well as  \cite{Wu21, Wu21a, Wu21b}),
\item[-] its close links with higher Auslander--Reiten theory
(to be illustrated below, cf. also \cite{Wu21, Wu21a}). 
\end{itemize} 

In these notes, after an informal illustration of the main notions and a
quick reminder on dg (=differential graded) algebras
and their derived categories, our first aim will be to present many examples of dg
algebras and morphisms endowed with Calabi--Yau structures, respectively relative
Calabi--Yau structures. These will be obtained using (relative) Calabi--Yau completions
and should illustrate the relevance of these for (higher) Auslander--Reiten theory.
Our second aim will be to sketch the foundations of the subject for which we
will need to recall the necessary material on Hochschild and cyclic homology.

\section*{Acknowledgments}
The authors thank the organizers of the ICRA 2020 for the excellent
job they have done. They are indebted to Yilin Wu for help with the
figures, the diagrams and the format.
They are grateful to Xiaofa Chen, Junyang Liu and an anonymous
referee for carefully reading previous versions of these notes and helping to 
improve their readability and to reduce the number of typos and mistakes. 

\section{Intuition and first examples}
\label{s:intuition}

In this purely introductory section, we informally discuss the key notions
and constructions to be developed in the sequel. 

There is a close analogy (which actually goes deeper) between
the notion of orientation of a (real, smooth) $n$-dimensional manifold $M$ and the notion of 
(absolute) $n$-Calabi--Yau structure on a dg algebra $A$. This extends to a 
relative setting\footnote{Notice that we do not write \emph{the} relative
setting because the setting considered here is \emph{not} the relative one
in the usual sense of algebraic geometry.} where we obtain a close analogy between an $n$-dimensional
manifold {\em with boundary} $\partial M \subset M$ both endowed with
compatible orientations and a {\em relative} $n$-Calabi--Yau structure
on a {\em morphism} $B \to A$ of dg algebras.

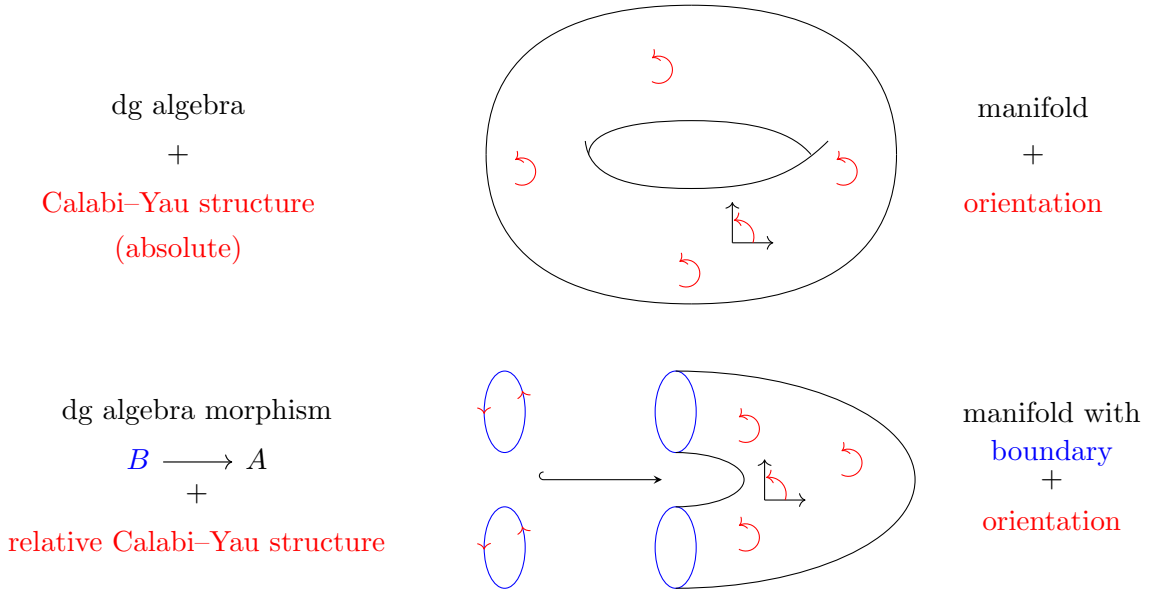
\begin{center}
	\begin{tikzpicture}[scale=0.9]
		\node (a) at (-7.5,0.7) {
			\text{dg algebra}
		};

		\node (b) at (-7.5,-0.7) {
			\color{red}\text{Calabi--Yau structure}
		};
		
		\node (c) at (-7.5,-1.4) {
			\color{red}\text{(absolute)}
		};
		\node at (-7.5,0) {+};
		
		\draw (-3,0) .. controls (-3,2) and (-1,2.2) .. (0,2.2);
		\draw[xscale=-1] (-3,0) .. controls (-3,2) and (-1,2.2) .. (0,2.2);
		\draw[rotate=180] (-3,0) .. controls (-3,2) and (-1,2.2) .. (0,2.2);
		\draw[yscale=-1] (-3,0) .. controls (-3,2) and (-1,2.2) .. (0,2.2);
		
		\draw (-1.55,.2) .. controls (-1.5,-0.3) and (-1,-0.5) .. (0,-.5) .. controls (1,-0.5) and (1.5,-0.3) .. (2,0.2);
		
		\draw (-1.5,0) .. controls (-1.5,0.3) and (-1,0.5) .. (0,.5) .. controls (1,0.5) and (1.5,0.3) .. (1.75,0);

		\draw[->,red] ([shift=(-45:0.6)](-3,0) arc (-120:120:0.2cm);
		
		\draw[->,red] ([shift=(-45:0.6)](1.7,0) arc (-120:120:0.2cm);
		
		\draw[->,red] ([shift=(-45:0.6)](-1,1.5) arc (-120:120:0.2cm);
		
		\draw[->,red] ([shift=(-45:0.6)](-0.6,-1.5) arc (-120:120:0.2cm);
		
		\node (d) at (5,0.7) {
			\text{manifold}
		};

		\node (e) at (5,-0.7) {
			\color{red}\text{orientation}
		};
		
		\node at (5,0) {+};
		
		\draw[->] (0.6,-1.3)--(1.2,-1.3);
		\draw[->] (0.6,-1.3)--(0.6,-0.7);
		
		\draw[->,red] ([shift=(0:0.6)]0.3,-1.3) arc (-20:100:0.25cm);
	\end{tikzpicture}
	
\end{center}

\vspace*{0.2cm}
\begin{center}
	\begin{tikzpicture}[scale=0.9]
        \node (a) at (-7.5,0) {
			\text{dg algebra morphism}
		};
	    \node (a) at (-7.5,-0.7) {
		\begin{tikzcd}
			\color{blue}B\arrow[r]&A
		\end{tikzcd}
	};
		\node at (-7.5,-1.2) {+};
		\node (b) at (-7.5,-1.9) {
			\color{red}\text{relative Calabi--Yau structure}
		};

		\draw[blue] (-3,0) ellipse (0.3 and 0.6 );
        \draw[->,red] (-3,0) ++ (30:0.3 and 0.6) arc (30:30.9:0.3 and 0.6);
        \draw[->,red] (-3,0) ++ (180:0.3 and 0.6) arc (180:183:0.3 and 0.6);
		
		\draw[blue] (-3,-2) ellipse (0.3 and 0.6 );
		\draw[->,red] (-3,-2) ++ (30:0.3 and 0.6) arc (30:30.9:0.3 and 0.6);
		\draw[->,red] (-3,-2) ++ (180:0.3 and 0.6) arc (180:183:0.3 and 0.6);

		\draw[blue] (-0.5,0) ellipse (0.3 and 0.6 );
		
		\draw[blue] (-0.5,-2) ellipse (0.3 and 0.6 );
		
		\draw (-0.5,-1) +(90:1cm and 0.4cm) arc (90:-90:1cm and 0.4cm);
		
		\draw (-0.5,-1) +(90:3.5cm and 1.6cm) arc (90:-90:3.5cm and 1.6cm);
		
		\node (d) at (5,0) {
			\text{manifold with}
		};
   	   \node at (5,-0.6) {
		\text{\color{blue}{boundary}}
 	   };
		
		\node (e) at (5,-1.6) {
			\color{red}\text{orientation}
		};
		
		\node at (5,-1) {+};
		
		\draw[->] (0.8,-1.3)--(1.4,-1.3);
		\draw[->] (0.8,-1.3)--(0.8,-0.7);
		
		\draw[->,red] ([shift=(0:0.6)]0.5,-1.3) arc (-20:100:0.25cm);
		
		\draw[->,red] ([shift=(-45:0.6)](0,0) arc (-120:120:0.2cm);
		
		\draw[->,red] ([shift=(-45:0.6)](0,-1.6) arc (-120:120:0.2cm);
		\draw[->,red] ([shift=(-45:0.6)](1.5,-0.5) arc (-120:120:0.2cm);
		
		\draw[right hook-stealth] (-2.5,-1)--(-0.7,-1);
	\end{tikzpicture}
	
\end{center}

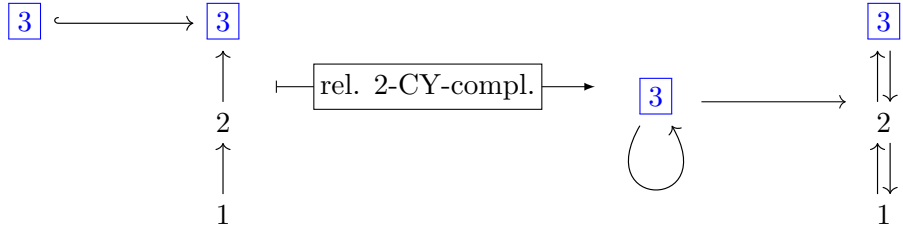
\begin{figure}

\begin{center}
\begin{tikzpicture}[scale=1]
	\node (a) at (0,0){
		\begin{tikzcd}
			\color{blue}\boxed{3}\arrow[rr,hook]&&\color{blue}\boxed{3}\\
			&&2\arrow[u]\\
			&&1\arrow[u]\\
		\end{tikzcd}
	};
\node (b) at (7,0){
\begin{tikzcd}
	\color{blue}\boxed{3}\arrow[out=240,in=300,loop,black]
\end{tikzcd}
};

\node(c) at (10,0){
\begin{tikzcd}
	\color{blue}\boxed{3}\arrow[d,shift left=0.5ex]\\
	2\arrow[d,shift left=0.5ex]\arrow[u,shift left=0.5ex]\\
	1\arrow[u,shift left=0.5ex]\\
\end{tikzcd}
};

\draw (2.5,1) rectangle (5.5,0.4) node[pos=.5] {rel. \textcolor{black}{2}-CY-compl.};
\draw[|-](2,0.7)--(2.5,0.7);
\draw[-latex](5.5,0.7)--(6.2,0.7);
\draw[->] (b)+(0.6,0.5)--(9.5,0.5);
\end{tikzpicture}
\end{center}

\caption{Relative $2$-Calabi--Yau completion}
\label{figure1}
\end{figure}

 In differential geometry,
we obtain examples of oriented manifolds with boundary by taking
the conormal bundle of the inclusion of a submanifold. In the homological
theory of dg algebras, we obtain examples of morphisms endowed with
relative Calabi--Yau structures by forming the $n$-Calabi--Yau completion $B \to A$ of a 
morphism $i_0: B_0 \to A_0$ of smooth dg algebras over a field $k$. For example,
let us take $B_0$ to be $k$, the algebra $A_0$ to be the path algebra
of the quiver $1 \to 2 \to 3$ and $i_0$ the morphism given by the inclusion
of the vertex $3$ into this quiver as in Figure~\ref{figure1}. If we apply the
relative $2$-Calabi--Yau completion to this morphism, we find the algebra
morphism from the polynomial algebra $k[t]$ to the Auslander algebra
of the truncated polynomial algebra $k[x]/(x^3)$ which takes the unique
object to the indecomposable projective $P=k[x]/(x^3)$ and the indeterminate $t$ 
to the multiplication by $x$. 

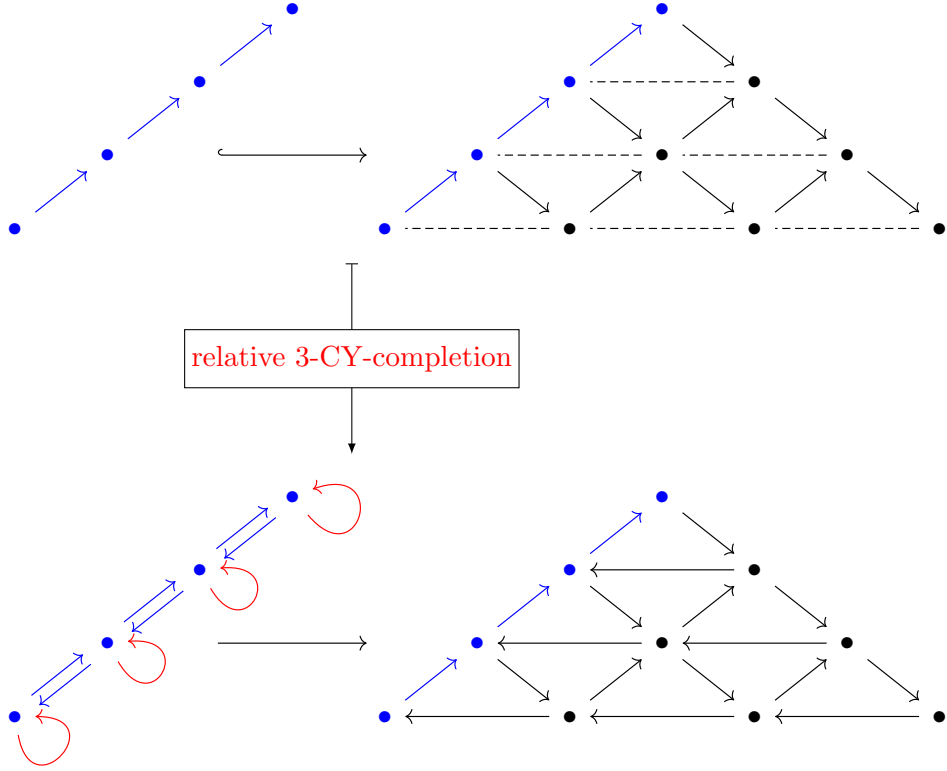
\begin{figure}
\begin{center}
		
\begin{tikzpicture}[scale=1.1]
	\node(a) at (0,0){
		\begin{tikzcd}[column sep=scriptsize, row sep=scriptsize]
			&&&\color{blue}\bullet&&&&\color{blue}\bullet\arrow[dr]\\
			&&\color{blue}\bullet\arrow[ur,blue]&&&&\color{blue}\bullet\arrow[ur,blue]\arrow[dr]&&\bullet\arrow[dr]\arrow[ll,dashed,no head]\\
			&\color{blue}\bullet\arrow[ur,blue]&\ \arrow[rr,hook]&&\ &\color{blue}\bullet\arrow[ur,blue]\arrow[dr]&&\bullet\arrow[dr]\arrow[ur]\arrow[ll,dashed,no head]&&\bullet\arrow[dr]\arrow[ll,dashed,no head]\\
			\color{blue}\bullet\arrow[ur,blue]&&&&\color{blue}\bullet\arrow[ur,blue]&&\bullet\arrow[ur]\arrow[ll,dashed,no head]&&\bullet\arrow[ur]\arrow[ll,dashed,no head]&&\bullet\arrow[ll,dashed,no head]\\
		\end{tikzcd}
		
	};
\node(b) at (0,-6){
	\begin{tikzcd}[column sep=scriptsize, row sep=scriptsize]
		&&&\color{blue}\bullet\arrow[dl,shift left=0.5ex,blue]\arrow[out=310,in=20,loop,red]&&&&\color{blue}\bullet\arrow[dr]\\
		&&\color{blue}\bullet\arrow[ur,shift left=0.5ex,blue]\arrow[dl,shift left=0.5ex,blue]\arrow[out=300,in=5,loop,red]&&&&\color{blue}\bullet\arrow[ur,blue]\arrow[dr]&&\bullet\arrow[dr]\arrow[ll]\\
		&\color{blue}\bullet\arrow[ur,shift left=0.5ex,blue]\arrow[dl,shift left=0.5ex,blue]\arrow[out=300,in=5,loop,red]&\ \arrow[rr]&&\ &\color{blue}\bullet\arrow[ur,blue]\arrow[dr]&&\bullet\arrow[dr]\arrow[ur]\arrow[ll]&&\bullet\arrow[dr]\arrow[ll]\\
		\color{blue}\bullet\arrow[ur,shift left=0.5ex,blue]\arrow[out=280,in=360,loop,red]&&&&\color{blue}\bullet\arrow[ur,blue]&&\bullet\arrow[ur]\arrow[ll]&&\bullet\arrow[ur]\arrow[ll]&&\bullet\arrow[ll]\\
	\end{tikzcd}
};
\draw (-3.5,-3) rectangle (0.5,-2.3) node[pos=.5] {\color{red}\text{relative 3-CY-completion}};
\draw[|-](-1.5,-1.5)--(-1.5,-2.3);
\draw[-latex](-1.5,-3)--(-1.5,-3.8);
\end{tikzpicture}
\end{center}
\caption{Example of a relative $3$-Calabi--Yau completion}
\label{figure2}
\end{figure}

In Figure~\ref{figure2}, we consider an example of a relative $3$-Calabi--Yau completion. The given morphism is the inclusion of the path algebra $A$ of
the linearly oriented $A_3$-quiver into its Auslander algebra. The $3$-Calabi--Yau
completion is a morphism from the $2$-dimensional Ginzburg algebra of
type $A_3$ to a relative Jacobian algebra, namely the one associated with
the ice quiver in the lower half of the figure endowed with the potential
given by the difference between the sum of the $3$-cycles rotating
clockwise and those rotating counter-clockwise.

\section{Complements on dg algebra resolutions}

Let $k$ be a field. Let $Q$ be a {\em graded quiver}, \ie a quiver $(Q_0, Q_1, s, t)$
endowed with a function $|?|: Q_1 \to \Z$. For $n\in \Z$ we write $Q^n_1$ for the
set of arrows $\alpha$ of $Q$ of degree $|\alpha|=n$. The associated path algebra $kQ$ 
then becomes a graded algebra if we define the degree of a path $a_1 \ldots a_l$ to
be the sum of the degrees of the composable arrows $a_i$. As a typical
example, let $A$ be the path algebra of the quiver
\[
\begin{tikzcd}
	1\arrow[r,"b"]&2\arrow[r,"a"]&3
\end{tikzcd}
\]
subject to the relation $ab=0$ (all arrows in degree $0$). Let $\tilde{A}$ be
the graded path algebra of the graded quiver
\[
\begin{tikzcd}
	1\arrow[r,"b"]\arrow[rr,"c",swap,bend right=35,red]&2\arrow[r,"a"]&3
\end{tikzcd}
\]
where $a$ and $b$ have degree $0$ and $c$ has degree $-1$. 
We endow the graded path algebra $\tilde{A}$ with the unique
algebra differential such that $d(c)=ab$. Then the algebra morphism
$\tilde{A} \to A$ taking $a$ to $a$, $b$ to $b$ and $c$ to zero is
compatible with the differential (that of $A$ being zero) and is
in fact a quasi-isomorphism of dg algebras. This is an example of
a cofibrant dg algebra resolution.

\begin{proposition} \label{prop: Cofibrant dg algebra resolution}
For each quotient $A=kQ/I$ of a path algebra with
finite $Q_0$, there is a dg algebra quasi-isomorphism
\[
\xymatrix{\tilde{A} \ar[r]^\eps & A}
\]
where $\tilde{A}=(k\tilde{Q},d)$ for a non positively graded quiver $\tilde{Q}$ with
$\tilde{Q}_0=Q_0$, $Q_1^0 = Q_1$ such that $\eps$ induces the identity
$\tilde{Q}_0=Q_0$ and the canonical projection $k\tilde{Q}_1 \to kQ_1$.
The morphism $\eps$ is a cofibrant dg algebra resolution.
\end{proposition}

In general, there is no `minimal choice' for $\tilde{A}$. However, there is
if $kQ/I$ is a monomial algebra, \cf \cite{Tamaroff21}, or if we work
in the setup of pseudocompact dg algebras, \cf \cite{VandenBergh15}.

\section{Derived categories and the inverse dualizing bimodule}

Let $k$ be a {\em perfect} field and $A$ a dg algebra. 
Recall that $\cc A$ denotes the
category of right dg $A$-modules, $\ch A$ the associated category up
to homotopy and $\cd A$ the derived category, \ie the localization of
$\cc A$ or $\ch A$ at the class of all quasi-isomorphisms. It is triangulated
with suspension functor $\Sigma: \cd A \to \cd A$ given by the shift
of dg modules and each short exact sequence of dg modules yields
a canonical triangle. The derived category $\cd A$ admits arbitrary
(set-indexed) coproducts. The perfect derived category $\per(A)$ is
its subcategory of compact objects. It equals the thick subcategory
of $\cd A$ generated by the free module $A_A$ of rank $1$. 

The {\em perfectly valued derived category $\pvd(A)$} is the full
subcategory of $\cd A$ whose objects are the dg modules $M$
whose underlying dg $k$-module is perfect. Equivalently, we may
require that we have
\[
\sum_{p\in\Z} \dim \Hom(P,\Si^p M) < \infty
\]
for each perfect dg module $P$ or that we have

\begin{equation*}
	\begin{split}
		\sum_{p\in \Z} \dim H^p(M) < \infty.
	\end{split}
\end{equation*}

The dg algebra $A$ is {\em proper} if it belongs to $\pvd(A)$ or, in other
words, if the sum of its homologies is finite-dimensional. For example,
if $A=A^0$ is concentrated in degree $0$, then $A$ is proper if
and only if $A^0$ is finite-dimensional and in this case, we have
\[
\pvd(A) = \cd^b(\mod A^0) \ko
\]
where $\mod A^0$ denotes the category of $k$-finite-dimensional
right $A^0$-modules.

The dg algebra $A$ is {\em connective} if $H^pA=0$ for all $p>0$.
In this case, the canonical morphism $\tau_{\leq 0} A \to A$ of dg
algebras is a quasi-isomorphism so that for most purposes, we
may assume that $A^p=0$ for $p>0$ when $A$ is connective.
If $A$ is connective, then $\cd A$ has a canonical $t$-structure
whose left (resp. right) aisle is formed by the dg modules $M$ such
that $H^p M=0$ for $p>0$ (resp. for $p<0$). Its heart identifies
with $\Mod H^0(A)$ via the restriction along the canonical
morphism $A \to H^0(A)$ (we may suppose $A^p=0$ for $p>0$). 
We say that $A$ is a {\em stalk dg algebra} or simply {\em stalk algebra}
if $H^p(A)=0$ for all $p\neq 0$. In this case, $A$ is linked to the dg algebra $H^0(A)$
concentrated in degree $0$ via the two dg algebra quasi-isomorphisms
\[
H^0(A) \la \tau_{\leq 0} A \to A.
\]

Let $A$ be a dg algebra. Its {\em enveloping algebra} is
\[
A^e = A \ten A^{op}
\]
so that right $A^e$-modules identify with $A$-$A$-bimodules via the rule
\[
m.(a \ten b) = (-1)^{|b||m|+|b||a|} bma.
\]
For a dg bimodule $M$, its {\em bimodule dual $M^\vee$} is defined by
\[
M^\vee = \RHom_{A^e}(M, A^e).
\]
We still view it as an object of $\cd(A^e)$ using the canonical isomorphism
$(A^e)^{op} = A^e$. More explicitly, the derived Hom-space
is computed using the outer bimodule structure on $A^e$ and the 
bimodule structure on $\RHom$ comes from the inner bimodule
structure on $A^e$. Notice that if $P$ is perfect in $\cd(A^e)$, we have
a canonical isomorphism
\[
P \iso (P^\vee)^\vee.
\]
Indeed, it suffices to check this for $P=A_A$ and then it is clear.

The dg algebra $A$ is (homologically) {\em smooth} if the {\em identity
bimodule} $\mbox{}_A A_A$ is perfect in $\cd(A^e)$. For example, if
$A=A^0$ is concentrated in degree $0$, then $A$ is smooth if and only
if $\mbox{}_A  A_A$ has a bounded resolution by finitely generated
projective bimodules. If $A=A^0$ is finite-dimensional, then $A$
is smooth if and only if it is of finite global dimension (here we use
the assumption that $k$ is perfect!). 

\begin{lemma}
\begin{itemize}
\item[a)] If $A$ is smooth, the category $\pvd(A)$ is Hom-finite.
\item[b)] If $Q$ is a finite non positively graded quiver, then
for any choice of algebra differential $d$ on the graded path
algebra $kQ$, the dg algebra $(kQ,d)$ is smooth.
\end{itemize}
\end{lemma}

Suppose that $A$ is a smooth dg algebra. Its {\em inverse dualizing
bimodule} is
\[
\Omega_A = A^\vee = \RHom_{A^e}(A, A^e).
\]

Denote by $D$ the duality $\Hom_k(?,k)$ over the ground field. The following
is Lemma~4.1 of \cite{Keller08d}
\begin{lemma}  \label{keylemma}
For $L\in \pvd(A)$ and $M\in \cd A$, we have a canonical isomorphism
\[
D\Hom_{\cd A}(L,M) \iso \Hom_{\cd A}(M\lten_A \Omega_A, L).
\]
\end{lemma}

\begin{corollary} The functor $S^{-1}=?\lten_A \Omega_A$ induces an
{\em inverse Serre functor} on $\pvd(A)$, \ie an autoequivalence such
that we have isomorphisms
\[
D\Hom_{\cd A}(L,M) \iso \Hom_{\cd A}(S^{-1} M,L)
\]
which are bifunctorial in $L,M\in \pvd(A)$.
\end{corollary}

\section{Calabi--Yau completions}

Fix a perfect field $k$. Let $A$ be a dg $k$-algebra. Fix an integer $n\in \Z$.
By definition, a {\em bimodule $n$-Calabi--Yau structure
on $A$} is an isomorphism
\[
\Sigma^n \Omega_A \iso A
\]
in $\cd(A^e)$.  Notice that by Lemma~\ref{keylemma}, the category
$\pvd(A)$ is then an $n$-Calabi--Yau category, \ie it is Hom-finite
and $\Sigma^n$ is a Serre functor. We say that $A$ is a {\em bimodule $n$-Calabi--Yau
dg algebra} if it is endowed with a bimodule $n$-CY structure.

\begin{exercise} Suppose that $A$ is bimodule $n$-Calabi--Yau and
$H^0(A)$ is not semi-simple.
\begin{itemize}
\item[a)] Show that if $A$ is connective and $H^0A$ is finite-dimensional,
then $H^p A\neq 0$ for infinitely many $p<0$.
\item[b)] Show that if $A$ is a stalk dg algebra, then $H^0 A$ is 
infinite-dimensional.
\end{itemize}
\end{exercise}

\subsection{Absolute Calabi--Yau completions}

Suppose that $B$ is a smooth dg algebra. Let $n\in \Z$ be an integer and
put
\[
\omega = \Sigma^{n-1} \Omega_B.
\]
We may and will assume that $\omega$ is cofibrant as a dg $B$-bimodule.
Following \cite{Keller11b}, we define the {\em $n$-Calabi--Yau completion}
of $B$ to be the tensor dg algebra
\[
\Pi_n B = T_B(\omega) = B \oplus \omega \oplus (\omega\ten_B \omega) \oplus
\cdots \oplus \omega^{\ten_B p} \oplus \cdots \ .
\]
It is not hard to check that up to quasi-isomorphism, $\Pi_n B$
does not depend on the choice of $\omega$ in its homotopy class.

\begin{theorem}[\cite{Keller11b, Yeung16, Keller18c}] The dg algebra
$\Pi_n B$ is smooth and carries a canonical bimodule $n$-Calabi--Yau structure.
In particular, the category $\pvd(\Pi_n B)$ is $n$-Calabi--Yau.
\end{theorem}

For example, consider $B=k$. Then $\Omega_B=k$ and $\Pi_n(B)=k[t]$, where
$t$ is of degree $1-n$ and $\Pi_n B$ carries the zero differential. In particular,
$\Pi_n k$ is concentrated in degree $0$ iff $n=1$. As another example,
let $Q$ be a connected non Dynkin quiver and $B=kQ$. Then the
$2$-Calabi--Yau completion $\Pi_2 B$ has its homology concentrated in
degree $0$ and is quasi-isomorphic to the preprojective algebra of $Q$.

Connective dg algebras form a particularly important class. It is natural
to ask to which extent this class is stable under forming CY-completions.
This question is not hard to answer: Suppose that $B$ is a smooth and
connective dg algebra. Define the {\em bimodule dimension} $d$ of
$B$ to be the supremum of the integers $p$ such that $H^p\Omega_B\neq 0$.
Then $\Pi_n B$ is connective if and only if $n\geq d+1$. 

Let us investigate $2$-CY-completions. Let $B$ be a smooth dg algebra.
For $n=2$, the bimodule $\omega$ is $\omega=\Sigma \Omega_B$ and the
functor
\[
?\lten_B \omega = \lten_B (\Sigma \Omega_B) = \Sigma \circ (?\lten_B \Omega_B)
\]
induces the composed functor
\[
\tau^{-1} = S^{-1} \Sigma
\]
in the perfect derived category $\pvd(B)$, where $S^{-1}$ is the inverse Serre
functor. We denote this functor by $\tau^{-1}$ because if $B$ is a
finite-dimensional algebra of finite global dimension, then it is
the inverse Auslander--Reiten translation of the category
$\pvd(B) = \cd^b(\mod B)$. For arbitrary smooth and proper $B$, 
the restriction of $\Pi_2(B)$ to a dg $B$-module is
\[
(\Pi_2(B)) |_B = T_B(\omega)|_B = B \oplus (B\lten_B \omega) \oplus 
(B\lten_B \omega\lten_B \omega) \oplus \cdots =
\bigoplus_{p\geq 0} \tau^{-p} B.
\]
Notice that in general, it will not be perfectly valued. Now suppose
that $B=kQ$ for a connected non Dynkin quiver $Q$. Then we know
that $\tau^{-p} B$ lies in $\mod B$ for all $p\geq 0$ so that
$\Pi_2(B)$ is a stalk dg algebra:
\[
\Pi_2(B) \iso H^0(\Pi_2 B) = T_B(H^0 \omega).
\]
It follows that $\Pi_2(B)$ is quasi-isomorphic to the classical
preprojective algebra of $Q$ by the description of this algebra
as a tensor algebra due to Baer--Geigle--Lenzing \cite{BaerGeigleLenzing87}.

If $Q$ is an arbitrary finite acyclic quiver and $B=kQ$, then $\Pi_2(B)$
can be described as the {$2$-dimensional Ginzburg algebra of $B$}.
For example, for the quiver 
\[
\begin{tikzcd}
	1\arrow[r,"b"]\arrow[rr,"c",swap,bend right=35,red]&2\arrow[r,"a"]&3
\end{tikzcd}
\]
the dg algebra $\Pi_2(B)$ is given by the graded quiver
\[
\begin{tikzcd}
	Q:&1\arrow[r,"b",shift left=0.5ex]\arrow[out=120,in=60,loop,red,"{t_{1}}"]&2\arrow[r,"a",shift left=0.5ex]\arrow[l,"{\overline{b}}",shift left=0.5ex]\arrow[out=120,in=60,loop,red,"{t_{2}}"]&3\arrow[l,"{\overline{a}}",shift left=0.5ex]\arrow[out=120,in=60,loop,red,"{t_{3}}"]\\
\end{tikzcd}
\]
with the arrows $a$, $b$, $\ol{a}$ and $\ol{b}$ in degree $0$
and the three loops $t_i$ in degree $-1$. The differentials of the
loops yield the preprojective relations:
\[
dt_1 = -\ol{b} b \ko d t_2 = b \ol{b} - \ol{a} a \ko d t_2 = a \ol{a}.
\]
Thus, we always have an isomorphism between $H^0(\Pi_2(B))$ and the
classical preprojective algebra but the homologies $H^p(\Pi_2(B))$
are non zero in infinitely many degrees $p<0$ unless all connected
components of $Q$ are non Dynkin. 

Let us now consider an example of a $3$-CY-completion. Consider
the Auslander algebra $B$ of the equioriented $A_3$-quiver given by
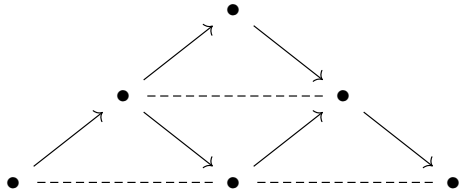
\[
\begin{tikzcd}
&&\bullet\arrow[dr]&&\\
&\bullet\arrow[ur]\arrow[dr]&&\bullet\arrow[dr]\arrow[ll,dashed,no head]&\\
\bullet\arrow[ur]&&\bullet\arrow[ur]\arrow[ll,dashed,no head]&&\bullet\arrow[ll,dashed,no head]
\end{tikzcd}
\]
Using a cofibrant dg algebra resolution of $B$
 (\cf  \ref{prop: Cofibrant dg algebra resolution}) it is not hard
to check that $\Pi_3(B)$ is quasi-isomorphic to the $3$-dimensional
Ginzburg algebra of the quiver with potential $(R,W)$, where $R$ is
the quiver of the `relation completion' of $B$
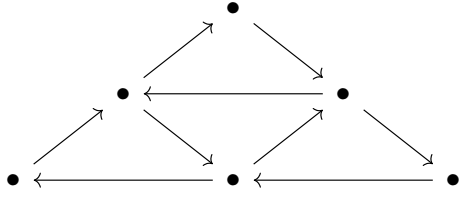
\[
\begin{tikzcd}
	&&\bullet\arrow[dr]&&\\
	&\bullet\arrow[ur]\arrow[dr]&&\bullet\arrow[dr]\arrow[ll]&\\
	\bullet\arrow[ur]&&\bullet\arrow[ur]\arrow[ll]&&\bullet\arrow[ll]
\end{tikzcd}
\]
and $W$ the difference of the sum of the $3$-cycles of $R$ turning
clockwise minus the unique $3$-cycle turning counterclockwise. For completeness,
let us recall how to construct the $3$-dimensional Ginzburg algebra
$\Gamma(R,W)$ associated with a quiver with potential $(R,W)$: starting
from $R$, construct a quiver $\tilde{R}$ as follows:
\begin{itemize}
\item[a)] for each arrow $a:i \to j$ of $R$, add an arrow $\ol{a}: j \to i$
of degree $-1$ and
\item[b)] for each vertex $i$ of $R$, add a loop $t_i: i\to i$ of degree $-2$.
\end{itemize}
Now define the differential on the graded path algebra $k\tilde{R}$ by
\[
d(t_i) = e_i \sum_{a\in R_1} (a \ol{a}-\ol{a} a) e_i
\]
for each vertex $i$ of $\tilde{R}$ and
\[
d\ol{a} = \partial_a W
\]
for each arrow $a$ of $R$. Here $\partial_a : HH_0(kR) \to kR$ is
the cyclic derivative defined on a path $p$ by
\[
\partial_a p = \sum_{p=uav} vu
\]
where the sum ranges over all decompositions $p=uav$ with 
paths $u$, $v$ of length $\geq 0$. For example, starting from
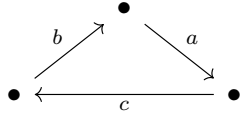
\[
\begin{tikzcd}
	&\bullet\arrow[dr,"a"]&\\
	\bullet\arrow[ur,"b"]&&\bullet\arrow[ll,"c"]
\end{tikzcd}
\]
with the potential $W=abc$ we obtain the quiver $\tilde{R}$
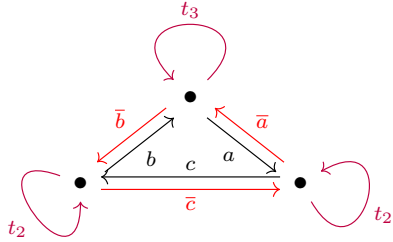
\[
\begin{tikzcd}
	&\bullet\arrow[dr,shift right=0.5ex,"a",swap]\arrow[dl,shift right=0.5ex,"{\overline{b}}",swap,red]\arrow[out=45,in=135,loop,purple,"{t_{3}}",swap]&\\
	\bullet\arrow[ur,shift right=0.5ex,"b",swap]\arrow[rr,shift right=0.5ex,"{\overline{c}}",swap,red]\arrow[out=160,in=270,loop,purple,"{t_{2}}",swap]&&\bullet\arrow[ll,shift right=0.5ex,"c",swap]\arrow[ul,shift right=0.5ex,"{\overline{a}}",swap,red]\arrow[out=300,in=30,loop,purple,"{t_{2}}",swap]
\end{tikzcd}
\]
with the differential determined by
\[
d(t_1)= c \ol{c} - \ol{b} b \ko
d(t_2)= a \ol{a} - \ol{c} c \ko
d(t_3)= b \ol{b} - \ol{a} a
\]
and
\[
d(\ol{a}) = bc \ko
d(\ol{b}) = ca \ko
d(\ol{c}) = ab.
\]

\section{Relative Calabi--Yau completions}

Let $k$ be a perfect field and $A$ and $B$ smooth dg $k$-algebras.
Let $f: B \to A$ be a dg algebra morphism. Recall that we do not
require $f$ to preserve the unit. A typical example would be the
inclusion of a finite-dimensional representation-finite algebra of finite
global dimension into its Auslander algebra. Let $n$ be an integer.

Following Yeung \cite{Yeung16}, we make the following definitions:
\begin{itemize}
\item[a)] the {\em relative inverse dualizing bimodule} is the bimodule dual of the cone
over the natural morphism
\[
A \lten_{B} A \to A
\]
considered as an object in $\cd(A^e)$. 
\item[b)] the $n$-dimensional {\em relative derived preprojective algebra} of $B \to A$
is
\[
\Pi_n(A,B) = T_A(\omega) \ko
\]
where $\omega$ is a cofibrant resolution of $\Si^{n-1} \Omega_{A,B}$.
\item[c)] the {\em relative $n$-Calabi--Yau completion} of $f: B \to A$ is the
canonical morphism of dg algebras
\[
\Pi_{n-1}(B) \to \Pi_n(A,B).
\]
\end{itemize}
Let us explain how the canonical morphism in c) is obtained: By construction,
we have a triangle in $\cd(A^e)$
\[
\Omega_{A,B} \to A^\vee \to (A \lten_B A)^\vee \to \Si\Omega_{A,B}.
\]
This yields morphisms
\[
\Omega_B = B^\vee \to A \lten_B B^\vee \lten_B A  \iso (A \lten_B A)^\vee \to 
\Si \Omega_{A,B}.
\]
The canonical morphism between the tensor algebras is induced by their composition.
We will define the notion of a relative (left) $n$-CY structure below. Assuming
it we can state the following theorem.

\begin{theorem}[Yeung \cite{Yeung16}, Bozec--Calaque--Scherotzke 
\cite{BozecCalaqueScherotzke20}] The dg algebra $\Pi_n(A,B)$ is
smooth and the morphism
\[
\Pi_{n-1}(B) \to \Pi_n(A,B)
\]
carries a canonical relative (left) $n$-CY structure.
\end{theorem}

Let us emphasize that all constructions and theorems generalize easily
from the setting of dg algebras and morphisms to that of dg categories and dg functors and
are proved in this generality in the references.

\section{Examples of relative $2$-Calabi--Yau completions}

\begin{figure}

\begin{center}
\begin{tikzpicture}[scale=1]
	\node (a) at (0,0){
		\begin{tikzcd}
			\color{blue}\boxed{3}\arrow[rr,hook]&&\color{blue}\boxed{3}\\
			&&2\arrow[u]\\
			&&1\arrow[u]\\
		\end{tikzcd}
	};
\node (b) at (7,0){
\begin{tikzcd}
	\color{blue}\boxed{3}\arrow[out=240,in=300,loop,blue]
\end{tikzcd}
};

\node(c) at (10,0){
\begin{tikzcd}
	\color{blue}\boxed{3}\arrow[d,shift left=0.5ex]\\
	2\arrow[d,shift left=0.5ex]\arrow[u,shift left=0.5ex]\\
	1\arrow[u,shift left=0.5ex]\\
\end{tikzcd}
};

\draw (2.5,1) rectangle (5.5,0.4) node[pos=.5] {rel. \textcolor{black}{2}-CY-compl.};
\draw[|-](2,0.7)--(2.5,0.7);
\draw[-latex](5.5,0.7)--(6.2,0.7);
\draw[->] (b)+(0.6,0.5)--(9.5,0.5);
\end{tikzpicture}
\end{center}

\caption{Relative $2$-Calabi--Yau completion}
\label{figure3}
\end{figure}
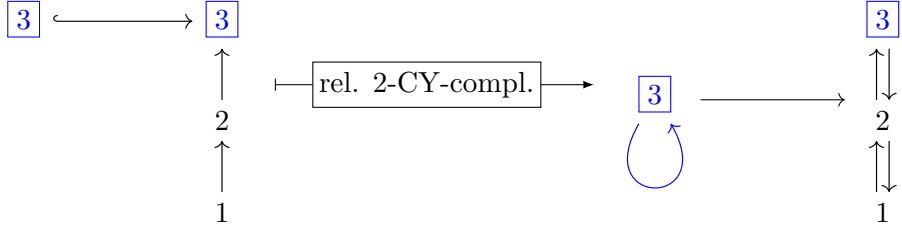

Let us recall the following example from section~\ref{s:intuition}:
Let us take $B$ to be $k$, the algebra $A$ to be the path algebra
of the quiver $1 \to 2 \to 3$ and $i$ the morphism given by the inclusion
of the vertex $3$ into this quiver as in Figure~\ref{figure3}. If we apply the
relative $2$-Calabi--Yau completion to this morphism, we find
a morphism of dg algebras $\Pi_1(B) \to \Pi_2(A,B)$ which, remarkably,
are both stalk algebras: the dg algebra $\Pi_1(B)$ is quasi-isomorphic
to the polynomial algebra $k[t]$ and the dg algebra $\Pi_2(A,B)$ to
the truncated polynomial algebra $k[x]/(x^3)$. The dg functor
$\Pi_1(B) \to \Pi_2(A,B)$ takes the unique
object to the indecomposable projective $P=k[x]/(x^3)$ and the indeterminate $t$ 
to the multiplication by $x$.

More generally, let us consider a finite quiver $Q$ and a subset
$F\subset Q_0$ of frozen vertices. We consider $F$ as a subquiver
with empty set of arrows. Then we have a natural algebra morphism
$kF \to kQ$. We then find that $\Pi_2(kQ,kF)$ is the {\em $2$-dimensional
relative Ginzburg algebra of $(Q,F)$}. For example, for the quiver
\[
\begin{tikzcd}
	\bullet\arrow[r,"b"]&\bullet\arrow[r,"a"]&\color{blue}\boxed\bullet
\end{tikzcd}
\]
we find that $\Pi_2(kQ,kF)$ is given by
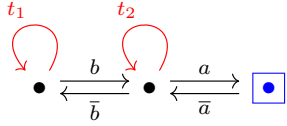
\[
\begin{tikzcd}
	\bullet\arrow[r,shift left=0.5ex,"b"]\arrow[out=60,in=130,loop,red,"{t_{1}}",swap]&\bullet\arrow[r,shift left=0.5ex,"a"]\arrow[l,shift left=0.5ex,"{\overline{b}}"]\arrow[out=60,in=130,loop,red,"{t_{2}}",swap]&\color{blue}\boxed\bullet\arrow[l,shift left=0.5ex,"{\overline{a}}"]
\end{tikzcd}
\]
where the arrows $a,b,\ol{a}$ and $\ol{b}$ are in degree $0$, the loops
$t_i$ in degree $-1$ and the differential is determined by
\[
d(t_1) = - \ol{b} b \ko d(t_2)=b \ol{b} - \ol{a} a.
\]
Thus, the algebra $H^0 \Pi_2(kQ, kF)$ is the {\em relative preprojective
algebra}: it has the same quiver as the classical preprojective algebra
but no relations at the frozen vertices. In our example here, we see
that remarkably, $\Pi_2(kQ,kF)$ is a {\em stalk algebra}. Here 
is another example:

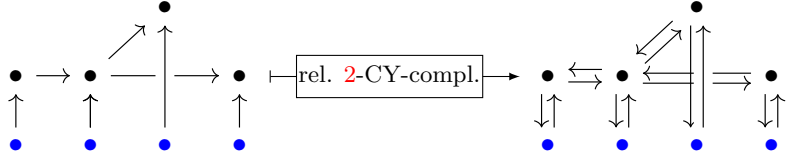
\begin{tikzpicture}[scale=0.7]
	\node (a) at (-1,0){
		\begin{tikzcd}[column sep=13, row sep=13]
			&&&\bullet\\
			&\bullet\arrow[r]&\bullet\arrow[rr]\arrow[ur]&&\bullet\\
			&\color{blue}\bullet\arrow[u]&\color{blue}\bullet\arrow[u]\arrow[u]&\color{blue}\bullet\arrow[uu,crossing over]&\color{blue}\bullet\arrow[u]
		\end{tikzcd}
	};
	
	\node (b) at (9,0){
		\begin{tikzcd}[column sep=13, row sep=13]
			&&&\bullet\arrow[dl,shift right=0.5ex]\arrow[dd,shift right=0.5ex]\\
			&\bullet\arrow[r,shift right=0.5ex]\arrow[d,shift right=0.5ex]&\bullet\arrow[rr,shift right=0.6ex]\arrow[d,shift right=0.5ex]\arrow[ur,shift right=0.5ex]\arrow[l,shift right=0.5ex]&&\bullet\arrow[ll,shift right=0.2ex]\arrow[d,shift right=0.5ex]\\
			&\color{blue}\bullet\arrow[u,shift right=0.5ex]&\color{blue}\bullet\arrow[u,shift right=0.5ex]&\color{blue}\bullet\arrow[uu,crossing over,shift right=0.5ex]&\color{blue}\bullet\arrow[u,shift right=0.5ex]
		\end{tikzcd}
	};
	
	\draw (2.5,-0.4) rectangle (6,0.4) node[pos=.5] {\footnotesize rel. \textcolor{red}{2}-CY-compl.};
	\draw[|-](a)--(2.5,0);
	\draw[-latex](6,0)--(6.7,0);
	
\end{tikzpicture}

Here we start from a `framed' quiver of type $D_4$ and obtain again a
relative preprojective algebra (in particular a stalk algebra) that could be called 
the `Nakajima algebra' because its (stable) representations (up to isomorphism)
with a given dimension vector  identify with the points of the corresponding
(regular) Nakajima quiver variety. 

As a final example, let us consider
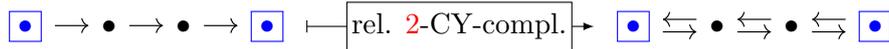
\begin{center}

\begin{tikzpicture}[scale=0.8]
\node (a) at (0,0){
\begin{tikzcd}[column sep=13, row sep=13]
	\color{blue}\boxed\bullet\arrow[r]&\bullet\arrow[r]&\bullet\arrow[r]&\color{blue}\boxed\bullet
\end{tikzcd}
};
\node (b) at (10,0){
\begin{tikzcd}[column sep=13, row sep=13]
	\color{blue}\boxed\bullet\arrow[r,shift right=0.5ex]&\bullet\arrow[r,shift right=0.5ex]\arrow[l,shift right=0.5ex]&\bullet\arrow[r,shift right=0.5ex]\arrow[l,shift right=0.5ex]&\color{blue}\boxed\bullet\arrow[l,shift right=0.5ex]
\end{tikzcd}
};

\draw (3.3,-0.4) rectangle (7,0.4) node[pos=.5] {rel. \textcolor{red}{2}-CY-compl.};
\draw[|-](a)--(3.3,0);
\draw[-latex](7,0)--(b);

\end{tikzpicture}	
\end{center}

Again, the dg algebra $\Pi_2(A,B)$ is a stalk algebra. Let us consider its
variant $\hat{\Pi}_2(A,B)$ where we replace the path algebra with the completed 
path algebra. This is easily seen to be isomorphic to the Auslander algebra
of the {\em Bass order}
\[
B_3 = \left[ \begin{array}{cc} R & Rx^3 \\ R & R \end{array} \right] \ko
\]
where $R=k[[x]]$ whose indecomposables are the $[R, Rx^i]$, $0\leq i\leq 4$.
For more information on Bass orders, we refer to page 72, of
\cite{Iyama01}.

\section{Examples of relative $3$-Calabi--Yau completions}

Consider the following example of a $3$-CY-completion:
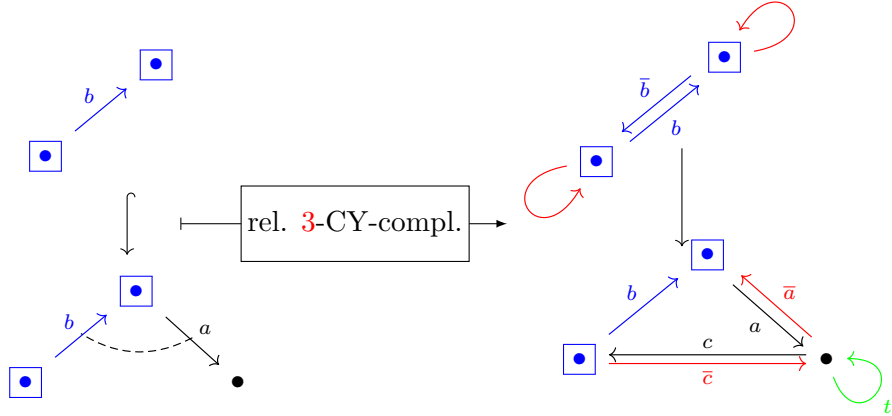
\begin{center}
\begin{tikzpicture}
\node(a) at (0,0){
\begin{tikzcd}[column sep=scriptsize, row sep=scriptsize]
	&\color{blue}\boxed\bullet\\
	\color{blue}\boxed\bullet\arrow[ur,blue,"b"]&&
\end{tikzcd}
};

\node(b) at (0,-3){
\begin{tikzcd}[column sep=scriptsize, row sep=scriptsize]
	&\color{blue}\boxed\bullet\arrow[dr,""{name=U},"a"]\\
	\color{blue}\boxed\bullet\arrow[ur,blue,""{name=V},"b"]&&\bullet\arrow[dashed,bend left=30,no head,from=U,to=V]
\end{tikzcd}
};

\node(c) at (7,0){
\begin{tikzcd}
	&\color{blue}\boxed\bullet\arrow[dl,blue,shift right=0.5ex,"\overline{b}",swap]\arrow[out=10,in=60,loop,red]\\
	\color{blue}\boxed\bullet\arrow[ur,blue,shift right=0.5ex,"b",swap]\arrow[out=190,in=240,loop,red]&&
\end{tikzcd}
};

\node(d) at (8,-3){
\begin{tikzcd}
	&\color{blue}\boxed\bullet\arrow[dr,shift right=0.5ex,"a",swap]\\
	\color{blue}\boxed\bullet\arrow[ur,blue,"b"]\arrow[rr,shift right=0.4ex,red,"\overline{c}",swap]&&\bullet\arrow[ll,shift right=0.3ex,"c",swap]\arrow[ul,shift right=0.5ex,red,"\overline{a}",swap]\arrow[out=290,in=0,loop,green,"t",swap]
\end{tikzcd}
};

\draw[right hook->](a)--(b);
\draw[->](7.3,-0.5)--(7.3,-1.8);
\draw (1.5,-2) rectangle (4.5,-1) node[pos=.5] {rel. \textcolor{red}{3}-CY-compl.};
\draw[|-](0.7,-1.5)--(1.5,-1.5);
\draw[-latex](4.5,-1.5)--(5,-1.5);

\end{tikzpicture}

\end{center}

Here, the morphism $B \to A$ on the left hand side is the embedding of
the path algebra of the quiver $A_2$  into its Auslander algebra.
On the right hand side, we have the (absolute) $2$-dimensional Ginzburg
algebra $\Pi_2(B)$ at the top and the $3$-dimensional relative
preprojective algebra $\Pi_3(A,B)$ at the bottom. It is quasi-isomorphic
to the relative $3$-dimensional Ginzburg algebra of the ice quiver
with potential $(R,W)$ with $R$ given by
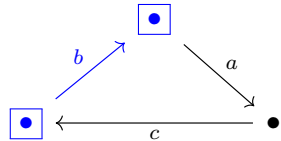
\[
\begin{tikzcd}
	&\color{blue}\boxed{\bullet}\arrow[dr,"a"]\\
	\color{blue}\boxed{\bullet}\arrow[ur,blue,"b"]&&\bullet\arrow[ll,"c"]
\end{tikzcd}
\]
and the potential $W=abc$. The construction of the relative $3$-dimensional
Ginzburg algebra is similar to that of the absolute one but
\begin{itemize}
\item[-] we do not add reversed arrows $\ol{b}$ for frozen arrows $b$ and
\item[-] we do not add loops $t_i$ for frozen vertices $i$.
\end{itemize}
In the example here, the relative $3$-dimensional algebra turns out
to be a stalk algebra, \ie quasi-isomorphic to the corresponding
relative Jacobian algebra, which is given by
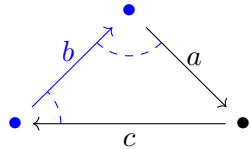
\[
\begin{tikzpicture}
	\node(a) at (0,0){
	$\color{blue}\bullet$};
\node(b) at (-1.5,-1.5){
	$\color{blue}\bullet$};
\node(c) at (1.5,-1.5){
	$\bullet$};
\draw[->](a)--(c) node [pos=0.6,above] {$ a $};
\draw[->](c)--(b) node [pos=0.5,below] {$ c $};
\draw[->,blue](b)--(a) node [pos=0.45,above] {$ b $};
\draw[dashed,blue] ([shift=(0:0.6)]-1.5,-1.5) arc (0:45:0.6cm);
\draw[dashed,blue] ([shift=(-45:0.6)]0,0) arc (-45:-135:0.6cm);
\end{tikzpicture}
\]
Note that the second example in section~\ref{s:intuition} is the analogous example
for equioriented $A_4$ instead of $A_2$. In the next section, we will exhibit
a general framework into which these examples fit.

\section{Higher Auslander algebras and the stalk property in dimension $\geq 3$}

Let $n\geq 1$ be an integer and $B$ a finite-dimensional algebra of
finite global dimension over a  perfect field $k$. Denote by $S$ the
Serre functor of the bounded derived category $\cd^b(B)=\cd^b(\mod B)$. Explicitly,
we have
\[
S= ?\lten_B DB.
\]
The composition of functors
\[
\xymatrix{
\mod B \ar[r]^{\mathrm{can}} &\cd^b(B) \ar[r]^{\Si^n S^{-1}} & \cd^b(B) \ar[r]^{H^0} &
\mod B}
\]
equals by definition $\tau_n^{-1}$, the {\em higher inverse Auslander-Reiten
translation} introduced by Iyama \cite{Iyama07a}. Define
\[
\cm=\add \{ \tau_n^{-p} B \; | \; p\geq 0\}.
\]
The motivating example is the case where $n=1$ and $B$ is the
path algebra $kQ$ of a Dynkin quiver $Q$. Then $B$ is $1$-dimensional
and $\mod(B)=\cm$ is $2$-dimensional and in the previous section,
we have seen examples where
$\Pi_3(\mod(B), \proj(B))$ is concentrated in degree $0$.

Put $\cb=\proj B$ and $\ca=\cm$. As shown by Iyama in \cite{Iyama07a},
if $\mod B$ admits an $n$-cluster-tilting object $M$, then the category
$\cm$ equals $\add(M)$. Let us assume that the category
$\ca$ (considered as a dg category concentrated in degree $0$) is
smooth (this assumption can be weakened to `local smoothness').
We consider the question of when the relative $(n+2)$-dimensional
derived preprojective algebra $\Pi_{n+2}(\ca,\cb)$ is a stalk category, \ie has its
homology concentrated in degree $0$. Recall that we have
\[
\Pi_{n+2}(\ca,\cb)=T_\cb(\omega) \ko
\]
where $\omega= \Si^{n+1} \Omega_{\ca,\cb}$ and $\Omega_{\ca,\cb}$ is a
cofibrant replacement of the bimodule dual of the cone over
\[
\ca \lten_\cb \ca \to \ca.
\]
For $M\in \cm$, let us abbreviate
\[
M^\wedge = \Hom_B(?,M)|_\ca
\]
considered as a finitely generated projective module over $\ca=\cm$.
We have the following key lemma.

\begin{keylemma} Suppose that $\cm$ is $n$-rigid, \ie we have
\[
\Ext^i_B(L,M)=0
\]
for all $L,M$ in $\cm$ and all $1\leq i\leq n-1$. Then for $M\in \cm$, we
have a canonical isomorphism
\[
M^\wedge \lten_\ca \omega \iso (\tau_n^- M)^\wedge.
\]
\end{keylemma}
A proof can be found in Proposition~8.6 of \cite{Wu21a}.

\begin{corollary} If $\cm$ is $n$-rigid, then the dg algebra
$\Pi_{n+2}(\cm, \proj B)$ is concentrated in degree $0$.
\end{corollary} 

Notice that the $n$-rigidity assumption holds if
\begin{itemize}
\item[-] $B$ is $n$-representation finite in the sense of
Iyama--Oppermann \cite{IyamaOppermann13} or
\item[-] $B^{op}$ is $n$-complete in the sense of \cite{Iyama11} or
\item[-] $B$ is $n$-representation infinite in the sense of
Herschend--Iyama--Oppermann \cite{HerschendIyamaOppermann14}
\end{itemize}

\section{Absolute and relative Calabi--Yau structures}

We have announced in the preceding sections that the $n$-CY-completion
$\Pi_n(B)$ of a smooth dg algebra carries a canonical (left)
$n$-Calabi--Yau structure and that the relative $n$-CY-completion
of a morphism $B \to A$ between smooth dg algebras carries
a canonical (left) relative $n$-CY-structure. Our aim in this
section is to define these structures and their right counterparts.

\subsection{A reminder on cyclic homology}
Let $k$ be a field and $A$ a $k$-algebra (associative, with $1$ but 
in general noncommutative). We consider $A$ as a bimodule over itself,
i.e.~as a module over the enveloping algebra $A^e=A\ten A^{op}$.
Recall that its {\em bar resolution $R$} is the complex whose component
in homological degree $p\geq 0$ is $A\ten A^{\ten p} \ten A$ and
whose differential $b'$ is defined for $p\geq 1$ by
\[
b'(a_0\ten a_1 \ten \cdots \ten a_p \ten a_{p+1}) = 
\sum_{i=0}^p (-1)^i a_0 \ten \cdots \ten a_i a_{i+1} \ten \cdots \ten a_{p+1}.
\]
The {\em Hochschild chain complex} is $CA=R\ten_{A^e} A$. Its
component in (homological) degree $p$ is isomorphic to $A^{\ten p}$ and
we write $b$ for its differential. It represents
the derived tensor product $A\lten_{A^e} A$ in the derived category of
vector spaces and its homology in (homological) degree $n$ is the
$n$th {\em Hochschild homology}
\[
HH_n(A) = H_n(R \ten_{A^e} A) = H_n(A \lten_{A^e} A)={\mathrm Tor}^{A^e}_n(A,A).
\]
For example, we find that we have $HH_0(A)=A/[A,A]$, where $[A,A]$ denotes
the subspace of $A$ generated by the commutators $[a,b]=ab-ba$, $a,b\in A$.

It is a remarkable fact due to Connes \cite{Connes83}
that the Hochschild chain complex admits a canonical `circle action
up to homotopy'. Algebraically, following Kassel \cite{Kassel85},
this is encoded by the structure of a {\em mixed complex},
i.e.~an action of the dg algebra
$\Lambda=k[\eps]/(\eps^2)$, where $\eps$ is of cohomological degree $-1$, endowed
with the differential $d=0$. Notice that $\Lambda$ is isomorphic to the singular
homology $H_*(S^1, k)$ of the circle $S^1$ and quasi-isomorphic to the dg
algebra of singular chains on $S^1$. It carries the structure of graded commutative
Hopf algebra with coproduct given by $\Delta(\eps) = 1 \ten \eps + \eps\ten 1$.
Connes has given explicit formulas for the $\Lambda$-action on $CA$ which
involve the unit $1\in A$. Since we need functoriality with respect to algebra
morphisms which do not necessarily preserve the unit, we replace $CA$
with a quasi-isomorphic complex $MA$, which we call the {\em mixed complex
of $A$}, and which is obtained as the total complex of the first two columns
of the Connes--Quillen bicomplex
\[
\begin{tikzcd}
\vdots\ar[d] & \vdots \ar[d] & \vdots \ar[d] \\
A^{\ten 2} \ar[d,"b"] & A^{\ten 2} \ar[d,"b'"]  \ar[l,"1-t"] & 
A^{\ten 2}  \ar[d,"b"] \ar[l, "N"]   & \ldots \ar[l,"1-t"] \\
A & A \ar[l,"1-t"] & A \ar[l,"N"] & \ldots\ .\ar[l, "1-t"] 
\end{tikzcd}.
\]
Here, as above, the symbol $b'$ denotes the differential of the bar resolution and
$b$ the differential of the Hochschild chain complex while $t$ is defined by
\[
t(a_1 \ten \cdots \ten a_n) = (-1)^{n-1} a_2 \ten \cdots \ten a_n \ten a_1
\quad\mbox{and}\quad
N=1 + t + \cdots + t^{n-1}.
\]
We write $d$ for the differential of $MA$ and $d': MA \to \Si^{-1} MA$ for
its endomorphism of degree $-1$ induced by $N$. Then we have
\[
d^2=0\ko {d'}^2=0 \quad\mbox{and}\quad d d' + d' d =0
\]
so that $MA$ becomes a dg module over $\Lambda$, where $\eps$ acts
via $d'$. In the sequel, we consider $MA$ as an object in the derived
category of mixed complexes, \ie the derived category $\cd\Lambda$ of
dg $\Lambda$-modules. It is then clear from the construction that
$MA$ is functorial in $A$ even with respect to algebra morphisms
which do not necessarily respect the unit, for example the inclusion
\[
\begin{tikzcd}
A \ar[r] & M_2(A) 
\end{tikzcd}
\]
sending $a\in A$ to the matrix
\[
 \left[\begin{array}{cc} a & 0 \\ 0 & 0 \end{array} \right].
 \]
 We have a canonical inclusion $CA \to MA$, which is a quasi-isomorphism
 because the second column is the augmented bar resolution, which is
 acyclic. Thus, we have an isomorphism
 \[
 A \lten_{A^e} A \iso MA
 \]
 in $\cd k$. In particular, we have isomorphisms
 \[
 HH_n(A) = H^{-n}MA = H^{-n}(MA \lten_\La \La) = H^{-n}(\RHom_\La(\La, MA)).
 \]
 One defines {\em cyclic homology} as the `coinvariants' of the homotopy circle action
 \[
 HC_n(A) = H^{-n}(MA \lten_\La k)
 \]
 and {\em negative cyclic homology} as the `invariants'
 \[
 HN_{n}(A) = H^{-n}(\RHom_\La(k, MA).
 \]
 Many authors write $HC^-_n(A)$ for $HN_n(A)$ but we prefer not
 putting the essential information into a barely visible exponent.
 The augmentation $\La \to k$ yields canonical maps
 \[
 \begin{tikzcd}
 HN_n(A) \ar[r] & HH_n(A) \ar[r] & HC_n(A).
 \end{tikzcd}
 \]
 For a morphism of algebras $B \to A$ (not necessarily preserving the unit), we define
 the {\em relative mixed complex $M(A,B)$} as the cone over the morphism $MB\to MA$.
 This yields the definitions of relative cyclic, Hochschild and negative homologies
 denoted respectively by $HC_n(A,B)$, $HH_n(A,B)$ and $HN_n(A,B)$.
 
 It is easy to extend these constructions from algebras to dg algebras and further
 to dg categories. For example, if $\ca$ is a dg category, the bar resolution of
 the identity bimodule $\ca(?, -)$ is given by the sum total complex of the double
 complex whose $p$th term is given by
 \[
 \coprod \ca(A_p, -) \ten \ca(A_{p-1},A_p) \ten \cdots \ten \ca(?,A_0),
 \]
 where the sum is over all sequences of objects $A_0$, \ldots\ , $A_p$ of $\ca$ and
 the differential is given by the same formula as in the case of an algebra.
 
 We recall the following facts from \cite{Keller98}: 
 
 \begin{theorem}
 \begin{itemize}
 \item[a)] Let $\ca$ and $\cb$ be dg categories and $F: \ca\to\cb$ a {\em Morita functor}, \ie
 a dg functor such that the restriction $F_*: \cd\cb\to\cd\ca$ is an equivalence.
 Then $MF: M\ca \to M\cb$ is an isomorphism in $\cd k$.
 \item[b)] Let
 \[
 \begin{tikzcd}
 0 \ar[r] & \ca \ar[r,"F"] & \cb \ar[r, "G"] & \cC \ar[r] & 0
 \end{tikzcd}
 \]
 be an {\em exact sequence} of dg categories, \ie the induced sequence
  \[
 \begin{tikzcd}
 0 \ar[r] & \cd\ca \ar[r,"F^*"] & \cd\cb \ar[r, "G^*"] & \cd\cC \ar[r] & 0
 \end{tikzcd}
 \]
 is an exact sequence of triangulated categories, where $F^*$ is left adjoint
 to the restriction $F_*$. Then there is a canonical triangle in $\cd\La$
 \[
 \begin{tikzcd}
 M\ca \ar[r,"MF"] & M\cb \ar[r, "MG"] & M\cC \ar[r] & \Si M\ca.
 \end{tikzcd}
 \]
 \end{itemize}
 \end{theorem}
 Notice that a) implies derived Morita invariance for the three invariants
 and that b) yields long exact sequences.

\subsection{Absolute left and right Calabi--Yau-structures}

Let $\ca$ be a smooth dg category. We have canonical isomorphisms
in $\cd k$
\[
\ca \lten_{\ca^e} {\ca} \iso \ca\lten_{\ca^e} (\ca^\vee)^\vee \iso 
\RHom_{\ca^e}(\ca^\vee, \ca).
\]
We deduce a canonical map
\[
HN_n(\ca) \to HH_n(\ca) \iso H^{-n}(\ca \lten_{\ca^e} \ca) \iso 
\Hom_{\cd \ca^e}(\ca^\vee, \Sigma^{-n} \ca).
\]

\begin{definition}[Kontsevich] A {\em left $n$-CY-structure on $\ca$} is a class
in $HN_n(\ca)$ whose image in $\Hom_{\cd \ca^e}(\ca^\vee, \Sigma^{-n} \ca)$
under the above map is an isomorphism.
\end{definition}

We see that if $\ca$ carries a left $n$-CY-structure, it is in particular
bimodule $n$-Calabi--Yau. The negative cyclic homology group
appears naturally when one studies the deformations of
bimodule $n$-CY-algebras, \cf \cite{ThanhofferVandenBergh12}.
One can show that each $n$-CY-completion of a smooth dg
category carries a canonical left $n$-CY-structure, \cf
\cite{Keller18c, BozecCalaqueScherotzke20}.

Now let $\ca$ be a proper dg category (\ie we have $\ca(X,Y)\in \per(k)$ for all $X,Y\in\ca$).
Let $D\ca^{op}$ denote the dg bimodule
\[
(X,Y) \mapsto D\ca(Y,X).
\]
We have canonical isomorphisms
\[
D(\ca\lten_{\ca^e}\ca) \iso \Hom_k(\ca\lten_{\ca^e} \ca, k) \iso
R\Hom_{\ca^e}(\ca, D\ca^{op}).
\]
We deduce morphisms
\begin{equation*}
	\begin{split}
		DHC_{-n}(\ca) \to D HH_{-n}(\ca) &\iso H^{-n}\RHom_{\ca^e}(\ca,D\ca^{op})\\
		&=
		\Hom_{\cd \ca^e}(\ca, \Si^{-n} D\ca^{op}).
	\end{split}
\end{equation*}

\begin{definition}[Kontsevich]  A {\em right $n$-CY-structure on $\ca$} is a class in
$DHC_{-n}(\ca)$ which yields an isomorphism $\ca \iso \Si^{-n} D \ca^{op}.$
\end{definition}

If $\ca$ carries a right $n$-CY-structure, we have the Serre duality formula
\[
D\ca(X,Y) \iso \Si^n \ca(Y,X)
\]
bifunctorially in $X,Y\in \ca$. In particular we see that $\ca$ is {\em perfectly
$n$-Calabi--Yau}, \ie the perfect derived category $\per(\ca)$ is $n$-Calabi--Yau as
a triangulated category. The definition above extends to `componentwise proper'
dg categories, \ie dg categories $\ca$ such that $H^p\ca(X,Y)$ is finite-dimensional
for all $p\in\Z$ and all $X,Y\in\ca$. These are important for our applications.
For example, if $Q$ is a finite acyclic quiver, then the canonical dg enhancement
\[
(\cc_Q)_{dg} = \per_{dg}(\Ga_{Q,0})/\pvd_{dg}(\Ga_{Q,0})
\]
of the classical cluster category $\cc_Q$ is componentwise proper
but not proper. Here the subscripts dg on the right hand side denote
the corresponding subcategories of the dg derived category
of $\Ga_{Q,0}=\Pi_3(kQ)$, \ie the canonical dg enhancement of the 
derived category.
For `componentwise proper' dg categories endowed with a
right $n$-CY-structure, the perfect derived category $\per(\ca)$ is still a 
Hom-finite $n$-CY triangulated category.

Suppose that $\ca$ is a dg category which is both smooth and proper.
Then, as shown by To\"en--Vaqui\'e in \cite{ToenVaquie07},
the Yoneda functor yields a quasi-equivalence
\[
\ca \iso \pvd_{dg}(\ca).
\]
Now suppose that moreover $\ca$ is augmented (in the sense of
\cite{Keller94}) and $\pvd_{dg}(\ca)$ 
is generated, as a triangulated category, by the image of
the restriction functor $\pvd(\ol{\ca}) \to \pvd(\ca)$, where
$\ca \to \ol{\ca}$ is the augmentation. Then we have a Morita
equivalence \cite{Keller06d}
\[
\ca^! \iso \pvd_{dg}(\ca).
\]
We deduce an isomorphism
\[
HN_n(\ca) \iso HN_n(\ca^!) \iso DHC_{-n}(\ca).
\]
\begin{proposition} This isomorphism establishes a bijection between
the left and the right $n$-Calabi--Yau structures on $\ca$.
\end{proposition}

Left CY-structures are inherited by dg localizations whereas right CY-structures
are inherited by full dg subcategories. More precisely, we have the
following proposition.

\begin{proposition}[\cite{Keller11b}] Let $F: \ca \to \cb$ be a dg functor. Suppose that
$F$ is a localization, \ie the functor $F^*: \cd\ca \to \cd\cb$ is a Verdier
localization. If $\ca$ is smooth, then so is $\cb$ and the image under
the induced morphism $HN_n(\ca) \to HN_n(\cb)$  of a left
$n$-CY-structure on $\ca$ is a left $n$-CY-structure on $\cb$.
\end{proposition}

\begin{remark} As a consequence, if $\Ga=\Ga_{Q,0}= \Pi_3(kQ)$ for an acyclic
quiver $Q$, then $\per_{dg}(\Ga)$ carries a left $3$-CY-structure
and so does its localization $(\cc_Q)_{dg}$, the dg cluster
category. This seems to be contradictory with the well-known
fact that the cluster category is $2$-Calabi--Yau as a
triangulated category. The explanation is that the $2$-Calabi--Yau
property of the triangulated category $\cc_Q$ comes from
a {\em right $2$-Calabi--Yau structure} on $(\cc_Q)_{dg}$ (\cf below).
The fact that $2\neq 3$ is not a contradiction because
$(\cc_Q)_{dg}$ is smooth (as a localization of $\per_{dg}(\Ga)$)
{\em but not proper} (only componentwise proper).
\end{remark}

\begin{theorem}[Brav--Dyckerhoff \cite{BravDyckerhoff19}] 
\label{thm: From left to right absolute CY-structures}
Let
$\ca$ be a smooth dg category. Each left $n$-CY structure on
$\ca$ yields a canonical right $n$-CY structure on $\pvd_{dg}(\ca)$.
\end{theorem}

For example, the canonical left $3$-CY structure on $\Ga$ as above
yields a right $3$-CY structure on $\pvd_{dg}(\Ga)$ which is responsible
for the $3$-CY property of the triangulated category $\pvd(\Ga)$. 
The right $3$-CY structure on $\pvd_{dg}(\Ga)$ yields the
right $2$-CY structure on $(\cc_Q)_{dg}$ via the connecting
morphism in cyclic homology associated with the exact sequence
\[
0 \to \pvd_{dg}(\Ga) \to \per_{dg}(\Ga) \to (\cc_Q)_{dg} \to 0 \ko
\]
\cf \cite{Keller98}.

\subsection{The derived category of morphisms}

Let $k$ be a perfect field and $A$ a dg $k$-algebra. Let $I=k\vec{A}_2$
(the letter $I$ stands for `interval'). Then the objects of $\cd(I^{op}\ten A)$
identify with morphisms $f:M_1 \to M_2$ of dg $A$-modules. Each object
gives rise to a triangle in $\cd(A)$
\[
\begin{tikzcd}
	M_{1}\arrow[r,"f"]&M_{2}\arrow[r]&C(f)\arrow[r]&\Si M_{1}
\end{tikzcd}
\]
functorial in $(M_1 \to M_2)$ considered as an object of
$\cd(I^{op}\ten A)$. For two objects $f: M_1 \to M_2$ and $f': M'_1 \to M'_2$,
consider a diagram whose rows are triangles of $\cd(A)$
\[
\begin{tikzcd}
	M_{1}\arrow[r,"f"]\arrow[d,dotted,"a",swap]&M_{2}\arrow[r]\arrow[d,"b",swap]&C(f)\arrow[r]\arrow[d,dotted,"c",swap]&\Si M_{1}\arrow[d,dotted,"\Si a",swap]\\
	M'_{1}\arrow[r,"f'"]&M'_{2}\arrow[r,"g'"]&C(f')\arrow[r]&\Si M'_{1}
\end{tikzcd}
\]
It is well-known (and easy to check) that for a given morphism $b$, there are
morphisms $a$ and $c$ making the diagram commutative if and only if
we have $g' b a=0$ and that in this case, the pair $(a,b)$ lifts to
a morphism of $\cd(I^{op}\ten A)$. This statement is refined by the following lemma

\begin{lemma} \label{lemma:morphisms}

We have a canonical isomorphism bifunctorial in $f,f' \in \cd(I^{op}\ten A)$
\[
\RHom_{I^{op}\ten A}(f,f') \iso\fib(\RHom_A(M_2, M_2') \to \RHom_A(M, C(f')) \ko
\]
where we write $\fib(g)$ for $\Sigma^{-1} C(g)$.
\end{lemma} 

A proof may be found in Lemma~3.1 of \cite{Wu21a}.

\subsection{Definition of relative Calabi--Yau structures}

Let $n\in\Z$ be an integer and $f: B \to A$ a morphism between smooth dg
algebras. With each class in relative Hochschild homology $HH_n(A,B)$, we
will associate a morphism of triangles of $\cd(A^e)$ 
\begin{equation} \label{eq: left CY-structure}
\tikzset{
	labl/.style={anchor=south, rotate=270, inner sep=.5mm}
}
\begin{tikzcd}[column sep=3.7em,row  sep=3em]
	\fib(\mu)\arrow[r,"\nu"]&A\lten_{B}A\arrow[r,"\mu"]&A\arrow[r]&\ldots\\
	\Si^{1-n}A^{\vee}\arrow[r,"\Si^{1-n}\mu^{\vee}"]\arrow[u,dotted]&
	\Si^{1-n}(A\lten_{B}A)^{\vee}\arrow[r]\arrow[u,dotted]&
	\Si^{1-n}\cof(\mu^{\vee})\arrow[r]\arrow[u,dotted]&\ldots
\end{tikzcd}
\end{equation}
For this, we observe that by Lemma~\ref{lemma:morphisms}, we have an isomorphism
\[
\RHom(\mu^\vee, \nu) \iso \fib(
\RHom_{A^e}((A\lten_B A)^\vee, A\lten_B A) \to \RHom_{A^e}(A^\vee, A)).
\]
Moreover, we have canonical isomorphisms
\[
A \lten_{A^e} A  \iso \RHom_{A^e}(A^\vee,A) 
\]
and
\[
A \lten_{B^e} A \iso (A\lten_B A) \lten_{A^e} (A \lten_B A) \iso
 \RHom_{A^e}((A\lten_B A)^\vee, A\lten_B A).
 \]
Thus, we get the following chain of morphisms
\begin{align*}
\RHom_{I^{op}\ten A^e}(\mu^\vee, \nu) & \iso  \fib(A \lten_{B^e} A  \to A \lten_{A^e} A) \\
                                    & \la  \fib(B \lten_{B^e} B \to A \lten_{A^e} A) \\
                                    & \iso \Si^{-1} HH(A,B) \\
                                    & \la \Si^{-1} HN(A,B) \ko
\end{align*}
where $HH(A,B)$ denotes the relative Hochschild complex of $B \to A$. 
By taking homology in homological degree $n$, we see that with each
class in $HN_n(A,B)$, there is associated a natural morphism of
triangles \ref{eq: left CY-structure}.

\begin{definition}[Brav--Dyckerhoff \cite{BravDyckerhoff19}] 
A {\em relative left $n$-CY-structure on $B\to A$} is a class in $HH_n(A,B)$ 
whose associated morphism of triangles \ref{eq: left CY-structure} is
an isomorphism.
\end{definition}

Notice that for $B=0$, we recover the absolute notion.
If $A$ and $B$ are concentrated in degree $0$, one easily deduces that
$A$ is bimodule internally $n$-Calabi--Yau (in the sense of Pressland \cite{Pressland17})
with respect to the image $e$ of $1_B$ in $A$. 

Suppose that $B \to A$ carries a left relative $n$-CY structure. Put $\ca=\pvd_{dg}(A)$
and $\cb=\pvd_{dg}(B)$ and let $r: \ca \to \cb$ be the restriction functor.

\begin{theorem}[Brav--Dyckerhoff \cite{BravDyckerhoff19}] The functor $r:\ca \to \cb$ inherits 
a canonical relative right $n$-CY structure, \ie there is a canonical class in 
$DHC_{1-n}(\cb,\ca)$ which yields an isomorphism of triangles
\[
\tikzset{
	labl/.style={anchor=south, rotate=270, inner sep=.5mm}
}
\begin{tikzcd}[column sep=3.7em,row  sep=3em]
	\Si^{n-1}\ca\arrow[r,"\Si^{n-1}r"]\arrow[d,"\sim"labl]&\Si^{n-1}\cb|_{\ca^{e}}\arrow[r]\arrow[d,"\sim"labl]&\Si^{n-1}\cof(r)\arrow[r]\arrow[d,"\sim"labl]&\ldots\\
	\fib(Dr)\arrow[r]&D\cb^{op}|_{\ca^{e}}\arrow[r,"Dr"]&D\ca\arrow[r]&\ldots
\end{tikzcd}
\]
where we write $\cof(r)$ instead of $C(r)$ to emphasize the duality with $\fib(Dr)$.
\end{theorem}

We refer to \cite{BravDyckerhoff19} for the construction of the class in $DHC_{1-n}(\cb,\ca)$
and that of the corresponding morphism of triangles in the theorem.
Concretely, for objects $L$ and $M$ of $\pvd(A)$, the above diagram becomes
\[
\tikzset{
	labl/.style={anchor=south, rotate=270, inner sep=.5mm}
}
\begin{tikzcd}[column sep=1.3em,row  sep=3em]
	&\RHom_{\ca}(L,\Si^{n-1}M)\arrow[r,"\text{res}"]\arrow[d,"\sim"labl]&\RHom_{\cb}(L,\Si^{n-1}M)\arrow[r]\arrow[d,"\sim"labl]&C(L,\Si^{n-1}M)\arrow[r]\arrow[d,"\sim"labl]&\ldots\\
	&DC(M,L)\arrow[r]&D\RHom_{\cb}(M,L)\arrow[r,"D\text{res}"]&D\RHom_{\ca}(M,L)\arrow[r]&\ldots
\end{tikzcd}
\]
Here, we write $C(M,L)$ for the cone over the morphism
\[
\RHom_\ca(M,L) \to \RHom_\cb(M,L)
\]
and $C(L, \Si^{n-1} M)$ for the cone over the morphism
\[
\RHom_\ca(L,\Si^{n-1}M) \to \RHom_\cb(L,\Si^{n-1}M).
\]
If we have $B=0$, we recover an isomorphism
\[
D\RHom_A(M,L) \iso \RHom_A(L, \Si^n M)
\]
so the theorem is a generalization of Theorem~\ref{thm: From left to right absolute CY-structures}.


\begin{thebibliography}{10}

\bibitem{BaerGeigleLenzing87}
Dagmar Baer, Werner Geigle, and Helmut Lenzing, \emph{The preprojective algebra
  of a tame hereditary {A}rtin algebra}, Comm. Algebra \textbf{15} (1987),
  no.~1-2, 425--457.

\bibitem{BozecCalaqueScherotzke20}
Tristan Bozec, Damien Calaque, and Sarah Scherotzke, \emph{Relative critical
  loci and quiver moduli}, arXiv:2006.01069 [math.AG].

\bibitem{BravDyckerhoff19}
\bysame, \emph{Relative {C}alabi-{Y}au structures}, Compos. Math. \textbf{155}
  (2019), no.~2, 372--412. \MR{3911626}
  

\bibitem{BravDyckerhoff18}
Christopher Brav and Tobias Dyckerhoff, \emph{Relative {C}alabi--{Y}au
  structures {II}: {S}hifted {L}agrangians in the moduli of objects},
  Selecta Math. (N.S.) \textbf{27} (2021), no.~4, article no.~63, 45 pp.
  arXiv:1812.11913 [math.AG].


\bibitem{Calabi57}
Eugenio Calabi, \emph{On {K}\"{a}hler manifolds with vanishing canonical
  class}, Algebraic geometry and topology. {A} symposium in honor of {S}.
  {L}efschetz, Princeton University Press, Princeton, N. J., 1957, pp.~78--89.

\bibitem{Connes83}
Alain Connes, \emph{Cohomologie cyclique et foncteurs {${\rm Ext}^n$}}, C. R.
  Acad. Sci. Paris S\'{e}r. I Math. \textbf{296} (1983), no.~23, 953--958.

\bibitem{GeissLeclercSchroeer11b}
Christof Gei\ss, Bernard Leclerc, and Jan Schr{\"o}er, \emph{Kac-{M}oody groups
  and cluster algebras}, Adv. Math. \textbf{228} (2011), no.~1, 329--433.

\bibitem{HerschendIyamaOppermann14}
Martin Herschend, Osamu Iyama, and Steffen Oppermann,
  \emph{{$n$}-representation infinite algebras}, Adv. Math. \textbf{252}
  (2014), 292--342. \MR{3144232}

\bibitem{Iyama01}
Osamu Iyama, \emph{Representation theory of orders}, Algebra---representation
  theory ({C}onstanta, 2000), NATO Sci. Ser. II Math. Phys. Chem., vol.~28,
  Kluwer Acad. Publ., Dordrecht, 2001, pp.~63--96.

\bibitem{Iyama07a}
\bysame, \emph{Higher-dimensional {A}uslander-{R}eiten theory on maximal
  orthogonal subcategories}, Adv. Math. \textbf{210} (2007), no.~1, 22--50.

\bibitem{Iyama11}
\bysame, \emph{Cluster tilting for higher {A}uslander algebras}, Adv. Math.
  \textbf{226} (2011), no.~1, 1--61.

\bibitem{IyamaOppermann13}
Osamu Iyama and Steffen Oppermann, \emph{Stable categories of higher
  preprojective algebras}, Adv. Math. \textbf{244} (2013), 23--68. \MR{3077865}

\bibitem{JensenKingSu16}
Bernt~Tore Jensen, Alastair~D. King, and Xiuping Su, \emph{A categorification
  of {G}rassmannian cluster algebras}, Proc. Lond. Math. Soc. (3) \textbf{113}
  (2016), no.~2, 185--212.

\bibitem{Kassel85}
Christian Kassel, \emph{Cyclic homology, comodules and mixed complexes}, J.
  Alg. \textbf{107} (1987), 195--216.


\bibitem{Keller94}
\bysame, \emph{Deriving {D}{G} categories}, Ann. Sci. {\'E}cole Norm. Sup. (4)
  \textbf{27} (1994), no.~1, 63--102.

\bibitem{Keller98}
\bysame, \emph{Invariance and localization for cyclic homology of {D}{G}
  algebras}, J. Pure Appl. Algebra \textbf{123} (1998), no.~1-3, 223--273.

\bibitem{Keller06d}
\bysame, \emph{On differential graded categories}, International Congress of
  Mathematicians. Vol. II, Eur. Math. Soc., Z\"urich, 2006, pp.~151--190.

\bibitem{Keller08d}
\bysame, \emph{Triangulated {C}alabi-{Y}au categories}, Trends in
  Representation Theory of Algebras (Zurich) (A.~Skowro\'nski, ed.), European
  Mathematical Society, 2008, pp.~467--489.

\bibitem{Keller11b}
\bysame, \emph{Deformed {C}alabi--{Y}au completions}, Journal f{\"u}r die reine
  und angewandte Mathematik (Crelles Journal) \textbf{654} (2011), 125--180,
  with an appendix by Michel~Van den Bergh.
  
\bibitem{Keller18c}
Bernhard Keller, \emph{Erratum to "{D}eformed {C}alabi--{Y}au completions"},
  arXiv:1809.01126 [math.RT, math.RA].
  

\bibitem{Leclerc16}
B.~Leclerc, \emph{Cluster structures on strata of flag varieties}, Adv. Math.
  \textbf{300} (2016), 190--228.

\bibitem{Pressland19}
Matthew Pressland, \emph{Calabi-{Y}au properties of {P}ostnikov diagrams},
Forum Math.~Sigma \textbf{10} (2022), article no.~e56, 31 pp.
  arXiv:1912.12475 [math.RT].

\bibitem{Pressland17b}
\bysame, \emph{A categorification of acyclic principal coefficient cluster
  algebras}, arXiv:1702.05352 [math.RT].

\bibitem{Pressland17}
\bysame, \emph{Internally {C}alabi-{Y}au algebras and cluster-tilting objects},
  Math. Z. \textbf{287} (2017), no.~1-2, 555--585.

\bibitem{Pressland20}
\bysame, \emph{Mutation of frozen {J}acobian algebras}, J. Algebra \textbf{546}
  (2020), 236--273. \MR{4033085}

\bibitem{Tamaroff21}
Pedro Tamaroff, \emph{Minimal models for monomial algebras}, Homology Homotopy
  Appl. \textbf{23} (2021), no.~1, 341--366. \MR{4185307}

\bibitem{Toen14}
Bertrand To\"{e}n, \emph{Derived algebraic geometry}, EMS Surv. Math. Sci.
  \textbf{1} (2014), no.~2, 153--240.

\bibitem{ToenVaquie07}
Bertrand To{\"e}n and Michel Vaqui{\'e}, \emph{Moduli of objects in
  dg-categories}, Ann. Sci. \'Ecole Norm. Sup. (4) \textbf{40} (2007), no.~3,
  387--444.

\bibitem{VandenBergh15}
Michel Van~den Bergh, \emph{Calabi-{Y}au algebras and superpotentials}, Selecta
  Math. (N.S.) \textbf{21} (2015), no.~2, 555--603.

\bibitem{ThanhofferVandenBergh12}
Michel Van~den Bergh and Louis de~Thanhoffer~de V\"olcsey, \emph{Calabi--{Y}au
  deformations and negative cyclic homology}, J.~Noncommut.~Geom.~\textbf{12} (2018),
  1255-1291. arXiv:1201.1520 [math.RA].

\bibitem{Wu21b}
Yilin Wu, \emph{Categorification of ice quiver mutation}, arXiv:2109.04503
  [math.RT].

\bibitem{Wu21}
\bysame, \emph{Relative {C}alabi--{Y}au structures in representation theory},
  Ph.~D.~thesis, Universit\'e Paris Cit\'e, December 2021.

\bibitem{Wu21a}
\bysame, \emph{Relative cluster categories and Higgs categories}, Advances \textbf{424},
1 July 2023, 109040. arXiv:2109.03707 [math.RT].

\bibitem{Yau78}
Shing~Tung Yau, \emph{On the {R}icci curvature of a compact {K}\"{a}hler
  manifold and the complex {M}onge-{A}mp\`ere equation. {I}}, Comm. Pure Appl.
  Math. \textbf{31} (1978), no.~3, 339--411.

\bibitem{Yeung16}
Wai-kit Yeung, \emph{Relative {C}alabi--{Y}au completions}, arXiv:1612.06352
  [math.RT].

\end{thebibliography}


\def\cprime{$'$} \def\cprime{$'$}
\providecommand{\bysame}{\leavevmode\hbox to3em{\hrulefill}\thinspace}
\providecommand{\MR}{\relax\ifhmode\unskip\space\fi MR }
\providecommand{\MRhref}[2]{%
  \href{http://www.ams.org/mathscinet-getitem?mr=#1}{#2}
}
\providecommand{\href}[2]{#2}

\end{document}